\documentclass[a4paper,fleqn]{cas-dc}

\usepackage[numbers]{natbib}

\usepackage{graphicx}
\usepackage{amsmath}
\usepackage{amsfonts}
\usepackage{amsthm}
\usepackage{cases}
\usepackage{wrapfig}
\usepackage{booktabs}
\usepackage{textcomp}
\usepackage{xcolor}
\usepackage{enumerate}
\usepackage{float}
\usepackage{caption}
\usepackage{lineno}

\def\tsc#1{\csdef{#1}{\textsc{\lowercase{#1}}\xspace}}
\tsc{WGM}
\tsc{QE}
\tsc{EP}
\tsc{PMS}
\tsc{BEC}
\tsc{DE}

\usepackage[mathscr]{euscript}

\begin{document}
\let\WriteBookmarks\relax
\def\floatpagepagefraction{1}
\def\textpagefraction{.001}

\shorttitle{Full Lyapunov Exponents spectrum with Deep Learning from single-variable time series}

\shortauthors{Carmen Mayora-Cebollero et~al.}

\title [mode = title]{Full Lyapunov Exponents spectrum with Deep Learning from single-variable time series}

\author[1]{Carmen Mayora-Cebollero}
\credit{cor1}

\author[1]{Ana Mayora-Cebollero}
\credit{-}

\author[2]{\'Alvaro Lozano}
\credit{-}

\author[1]{Roberto Barrio}
\credit{cor1}

\affiliation[1]{organization={IUMA, CoDy and Department of Applied Mathematics, Universidad de Zaragoza},
    city={Zaragoza},
    country={Spain}}

\affiliation[2]{organization={IUMA, CoDy and Department of Mathematics, Universidad de Zaragoza},
    city={Zaragoza},
    country={Spain}}

\cortext[cor1]{Corresponding author}

\begin{abstract}
    In this article we study if a  Deep Learning technique can be used to obtain an approximated value of the Lyapunov exponents
of a dynamical system. Moreover, we want to know if Machine Learning techniques are able, once trained, to provide the complete Lyapunov exponents spectrum
with just single-variable time series.  We train a Convolutional Neural Network and we use the resulting network to approximate the complete spectrum using the time series of just one variable from the studied systems (Lorenz system and coupled Lorenz system). The results are quite stunning as all the values are well approximated with only partial data. This strategy permits to speed up the complete analysis of the systems and also to study the hyperchaotic dynamics in the coupled Lorenz system.
\end{abstract}


\maketitle


\section{Introduction}
\label{sec:0}

Lyapunov Exponents (LEs) are a classical tool to study the behaviour of a dynamical system. For example, a positive maximum LE (MLE) means chaotic behavior, or hyperchaotic if the second LE is also positive. Moreover, with the second LE, some bifurcations, like period doubling bifurcations, can be inferred. One of the standard algorithms for the computation of LEs can be found in~\cite{WSSV85}, where the whole variables and the corresponding variationals are used. Other techniques~\cite{WSSV85,ROSENSTEIN1993117} use only one of the system variables time series, but only the maximum LE is obtained. Deep Learning (DL) techniques have been used to predict the LE of one-dimensional discrete dynamical systems directly~\cite{makarenko2018deep}  or as a complementary tool. In this last approach, DL is used for time series forecasting or to obtain, via data assimilation, a conjectured dynamical system, and later the Lyapunov spectrum is computed via classical methods \cite{golovko2003estimation, savitsky2015technique, dmitrieva2016method, bompas2020accuracy, pathak2017using}.

Deep Learning~\cite{HH19,G19} includes all the Machine Learning techniques that allow Deep Artificial Neural Networks (architectures built with layers of artificial neurons) to learn from data with several levels of abstraction. The activation of each neuron of a DL architecture (that is, its value) is computed applying a non-linear activation function to a linear combination of its inputs using some weights and a bias (that are known as the network trainable parameters). To fit all these parameters (in order to obtain the desired output for each input) a minimization problem is solved (a loss function is minimized respect to these parameters during a process known as training that uses training data). Unseen data, known as test data, is used to check that the network has learnt correctly and it is able to generalize to new samples.

The calculation of biparametric analysis with classical methods has allowed to study in detail the global dynamics of numerous systems~\cite{BS07,BS09,GG14}. Any improvement that can permit to faster and carry out more detailed studies could be useful. In the past few years, some authors have proposed to use Deep Learning as a new technique to analyse the behaviour of a dynamical system~\cite{makarenko2018deep,golovko2003estimation, savitsky2015technique, dmitrieva2016method, bompas2020accuracy, pathak2017using,ChaosDetCody23,BDNS20,LF20,CGRV22}. Such new technique can speed up these parametric studies.

In this paper, we apply DL to directly approximate the values of the LEs of two dynamical systems chosen as test examples: the classical Lorenz system~\cite{L63} and a coupled Lorenz system~\cite{GRASSI2009284}. Among all possible DL architectures, we have chosen the Convolutional Neural Network (CNN)~\cite{LBBH98}. Moreover, we want to show that, with just a short time series of one system variable, DL can be used to approximate all the LEs of a system. The notation ${\rm{LE}}_i$, where $i$ represents the ordered number of the exponent, will be used to indicate the exponent we are referring to (notice that ${\rm{LE}}_1$ corresponds to the MLE).

\paragraph{Lorenz system}

The Lorenz system~\cite{L63} is a classical three-dimensional continuous dynamical system given by the system of equations
\begin{equation}
	\left\{
	\begin{array}{rcl}
		\dot{x} & = & \sigma(y-x),\\
		\dot{y} & = & -xz + rx - y,\\
		\dot{z} & = & xy - bz,
	\end{array}
	\right.
	\label{LS_eq}
\end{equation}
where $(x,y,z)$ are the system variables and $(\sigma, r, b)$ are the bifurcation parameters ($\sigma$ is the Prandtl number, $r$ is the relative Rayleigh number, and $b$ is a positive constant). This is one of the seminal systems of chaotic dynamics. In the literature there are numerous papers~\cite{BS07,BS09} where its behaviour is described using classical LEs algorithms.

\paragraph{Coupled Lorenz system}
In order to increase the dimensionality of the test problem, we use two coupled Lorenz systems like in \cite{GRASSI2009284}:
\begin{equation}
	\left\{
	\begin{array}{rcl}
		\dot{x}_1 & = & \sigma(y_1-x_1),\\
		\dot{y}_1 & = & -x_1z_1 + r_1x_1 - y_1 + \lambda_1(x_2 - y_2),\\
		\dot{z}_1 & = & x_1y_1 - bz_1,\\
		\dot{x}_2 & = & \sigma(y_2-x_2),\\
		\dot{y}_2 & = & -x_2z_2 + r_2x_2 - y_2 + \lambda_2(x_1 - y_1),\\
		\dot{z}_2 & = & x_2y_2 - bz_2,\\
	\end{array}
	\right.
	\label{eqcls}
\end{equation}
where $(x_1,y_1,z_1,x_2,y_2,z_2)$ are the system variables, $(\sigma,$ $r_1,$ $r_2, b)$ are the bifurcation parameters, and $(\lambda_1, \lambda_2)$ are the coupling parameters. In what follows we consider $r_1=r$ and $r_2 = r-10$. In~\cite{GRASSI2009284} a dynamical study of the attractors of such system in the case of coupled equal Lorenz systems ($r_1=r_2$) can be found.

\medspace

This paper is organized as follows. In Section~\ref{sec:DL}, we introduce briefly the CNN architecture.  In Section~\ref{sec:lorenz}, we focus on the LEs prediction task for the Lorenz system. In Section~\ref{sec:coupledlorenz}, we show the results obtained in the prediction of LEs in the Coupled Lorenz system. Finally, in Section \ref{sec:Conclusions} we draw some conclusions.

All the DL experiments in this work have been performed with PyTorch~\cite{Py19}. The code has been executed on a Linux box with dual Xeon ES2697 with 128Gb of DDR4-2133 memory with a RTX2080Ti GPU.

\section{DL Techniques for Lyapunov Exponents Approximation}\label{sec:DL}

Deep Learning~\cite{HH19,G19} is the branch of Machine Learning that uses Deep Artificial Neural Networks  to learn from data with several levels of abstraction. These Artificial Neural Networks (ANNs) are formed by artificial neurons loosely inspired by their  biological counterparts that are organized in layers.

In a previous paper~\cite{ChaosDetCody23} we focused on the detection of chaotic behaviour in a dynamical system, but now we also want to quantify it and to be able to approximate all the Lyapunov Exponents using single-variable time series. In \cite{ChaosDetCody23} we used three ANN technologies:  the Multi-Layer Perceptron (MLP), Convolutional Neural Networks (CNNs) and Long Short-Term Memory networks (LSTMs). Here, we just use the CNN as it seems to work properly for this task. In the CNN design, we use a not very complicated structure and we do not perform hyperparameter optimization as we are mainly interested in studying the reliability of the methodology.

\begin{figure}[htb]
    \includegraphics[width=0.45\textwidth]{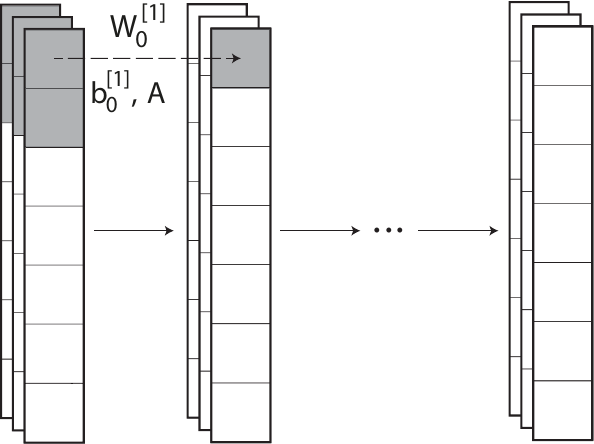}
    \caption{Simple graphic representation of the architecture of a 1D CNN with three channels in the depicted convolutional layers.}
    \label{ANNs}
\end{figure}

CNNs were originally developed for image recognition~\cite{LBBH98}, and are organized into convolutional and pooling layers that capture features and reduce dimensions, respectively. One of the main features of CNNs is that they share weights across multiple neurons~\cite{BBCV21}. Moreover, they can handle different input formats (vectors, matrices,...) depending on the type of convolution used. In this paper, the input data is in vector form, and therefore we use 1D CNNs, as depicted in Figure~\ref{ANNs}.

To exemplify how convolutional layers work, we show how to compute the value of the gray neuron in the second layer of the CNN in Figure~\ref{ANNs}, which is given by
\[
	x^{[1]}_{0,0} = \mathcal{A}\left(b^{[1]}_{0} + \sum_{j=0}^{1}\sum_{k=0}^{2}w_{j,k,0}^{[1]}\, x^{[0]}_{j,k}\right),
\]
where $x^{[l]}_{j,k}$ is the activation of neuron $j$ of channel $k$ at layer $l$ (first index for neurons, channels and layers is $0$), $W_{0}^{[1]} = (w^{[1]}_{j,k,0})_{j=0,1;\, k=0,2}$ is the weight matrix, $b^{[1]}_{0}$ is the bias vector, and $\mathcal{A}$ is the non-linear activation function.

In Section~\ref{sec:lorenz}, we detail the used CNN (number of convolutional layers, kernel size,\dots).

\section{LEs Approximation for the Lorenz System}\label{sec:lorenz}

In this section, two approaches are presented to approximate LEs with DL in the Lorenz system. On the first approach, the CNN network is trained (and validated) using information from few $r$-parametric lines, and its performance is tested. In the second approach, the same network architecture is trained (and validated) from scratch using time series with random parametric values from an $(r,b)$-parametric plane. These two approaches will show that DL techniques can be used to expand a partial classical study of the system (first approach) or to do the analysis from random data (second approach). The advantages and disadvantages of each approach are also compared.

In what follows, the numerical results of the performance of the DL process will be given as {mean\,$\pm$\,standard deviation} for $10$ randomly initialized networks.

\subsection{Non-Random Data}\label{sec:LS_NonRandom}

In our first test we consider uniparametric lines to create the datasets used to train the network for LEs approximation in the Lorenz system.

For this approach, the available data belongs to four $r$-parametric lines with $b\in\{2,2.4,8/3,2.8\}$ ($\sigma = 10$ for all the lines). For each line we consider $6,000$ different values of $r\in(0,300]$, which makes a total of $24,000$ time series. Such samples are obtained using the DOPRI5 integrator (a well known RK of order 5), but only the $x$-time series and the full LEs spectrum will be used. In order to obtain more precision, a transient integration is performed until time $t = 100,000$ with time step 0.01, later the integration is continued for $10,001$ time units with time step $0.001$. LEs used as the ground truth in the DL process are computed during this last integration. The time series are built with $1$ out of every $100$ of the last $100,000$ computed points. Similar samples are deleted according to the rule that two time series $p_{1}$ and $p_{2}$ are similar if $\Vert p_1 - p_2\Vert_{\infty}<10^{-4}$, and the remaining ones are normalized ($x$-coordinate is normalized linearly mapping its range to the interval $[0, 1]$; if the time series is constant over time, a random value between $0$ and $1$ is assigned). From the samples satisfying $b\in\{2,8/3\}$, $8,000$ are chosen randomly for the training set (batch size $128$). From the data on the $r$-parametric line with $b=2.4$, we select $2,000$ random points for validation (batch size $100$). Finally, $2,000$ random samples from the set $b=2.8$ are used for the test set (batch size $100$). Notice that, during training process only data from three lines is used: two lines corresponds to training (see light green lines in top-left panel of Figure~\ref{2DNonRandom_LS}) and another one to validation (see dark green line in the aforementioned figure).

The CNN architecture used for the task is inspired by that in~\cite{ChaosDetCody23} used for a chaos detection analysis (a classification task from the DL point of view). The network has only one input channel in the input layer as only the $x$-variable of the Lorenz system is used as known information. It has two convolutional layers: the first one with $15$ channels, kernel size $10$, stride $1$ and dilation equal to $2$; and the second one with $30$ channels, kernel size $5$, stride $1$ and dilation $4$. The ReLU activation function is applied after a bias term is added on each convolutional layer. Zero-padding and cropping are used to ensure that the length of the input sequences remains along the convolutional layers since the stride is $1$. A global average pooling layer is applied after the last convolutional layer and a readout layer with three networks (each one for the value of a LE) and bias term is stacked at the end. No activation function is applied in the output layer. Adam with learning rate $0.008$ is used as the optimizer algorithm. The ${\rm{L}}^2$-regularization with weight decay $10^{-5}$ is used to prevent overfitting. We use the Huber loss function~\cite{huber1964robust, hastie2009elements} given by
\[
\text{Huber Loss} = \dfrac{1}{N}\sum_{j=1}^{N}\,l_{j},\]
\[\text{ with } l_{j} = \left\{\begin{array}{rl}0.5\,(x_j-y_j)^2 & \text{ if }\vert x_{j}-y_{j}\vert<\delta,\\\delta\,(\vert x_{j}-y_{j}\vert-0.5\,\delta) & \text{otherwise},\end{array}\right.
\]
where $N$ is the batch size (for example, for our training dataset $N=128$) and we set $\delta=0.6$. Huber loss is less sensitive to outliers than the Mean Squared Error (MSE) loss (usually used in prediction tasks). In Huber loss the advantages of MSE loss and $L^{1}$-loss are combined. We expect that this fact will allow the CNN to focus more on fitting values near zero than those of large values (notice that a large error predicting a zero LE can produce a wrong detection of its behaviour). The number of epochs is $2,000$. An early stopping technique~\cite{Goodfellow-et-al-2016} is applied, so the used trained network is that with the values of weights and biases that give the lowest Huber loss value for the validation dataset while training.

The resulting value of the Huber loss for training dataset is $0.049\pm0.009$ (remember that the numerical results of the performance of the DL technique are given as mean $\pm$ standard deviation for $10$ randomly initialized networks). The corresponding values for validation and test datasets are $0.118\pm0.004$ and $0.113\pm0.012$, respectively. All mean values are close to zero and the standard deviations are small. The value of the loss function in the test set is a bit larger than that of the training set, however, we consider that there is not a large enough difference to discard the network because of overfitting. The data used for the analyses with the trained CNN are just time series of $x$-variable of length $1,000$.

In Subsection~\ref{ss1}, a $1$D analysis is performed to show that the trained network is able to generalize to an $r$-parametric line where it was not trained. In Subsection~\ref{ss2}, the analysis is extended to an $(r,b)$-biparametric plane.

\subsubsection{1D Analysis}
\label{ss1}
\begin{figure}[htb]
\centerline{\includegraphics[width=0.5\textwidth]{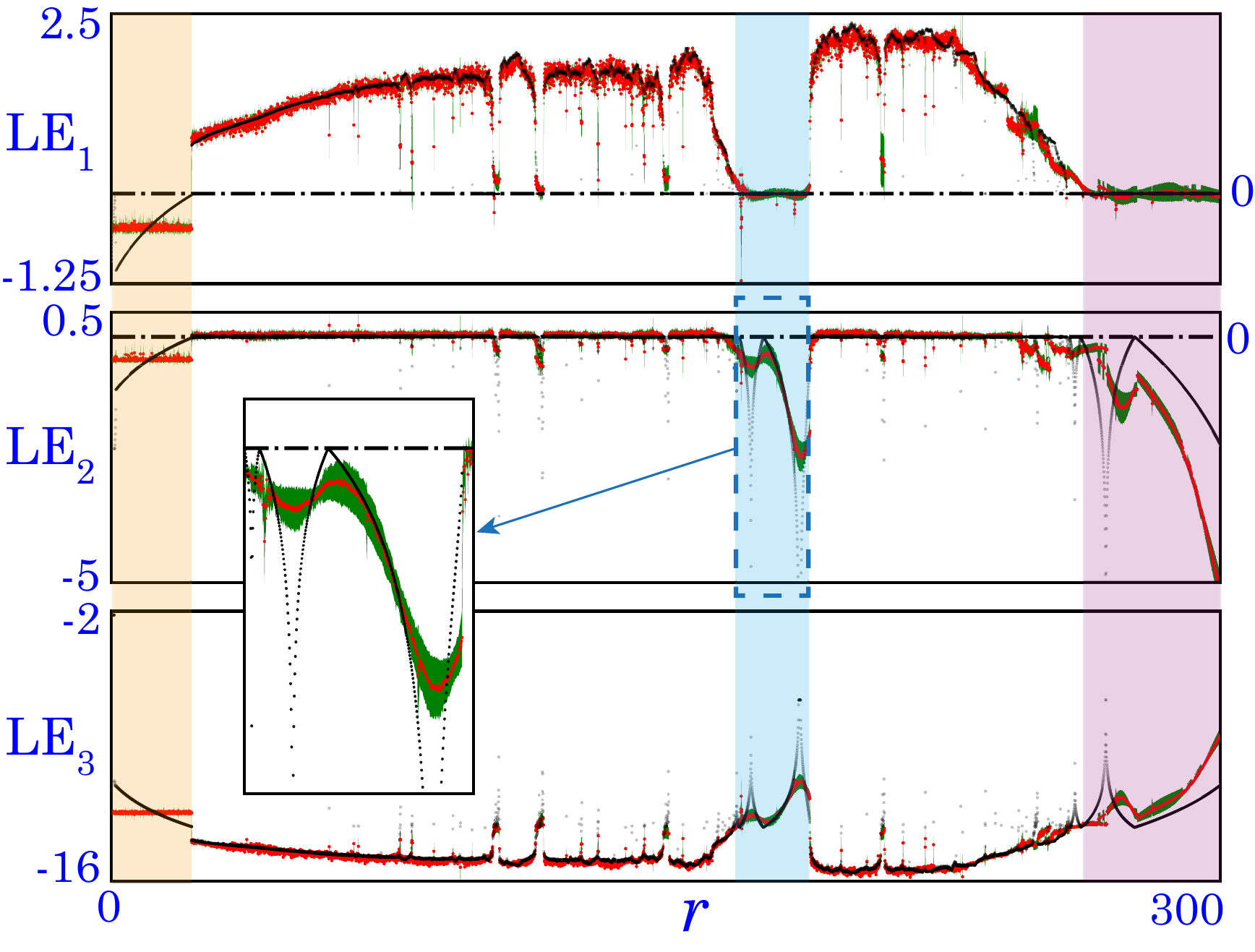}}
	\caption{1D analysis ($\sigma=10$, $b=2.2$) of Lyapunov Exponents in the Lorenz system when training the CNN with non-random data (Huber loss value $0.091\pm0.005$). Parts with orange and purple back colors correspond to the regions where the DL technique seems to give worst results. Region shaded in blue is used in a comparison with a subsequent analysis in Figure~\ref{1DRandom_LS}. (See the text for more details.)}
	\label{1DNonRandom_LS}
\end{figure}

In Figure~\ref{1DNonRandom_LS} the trained CNN has been used for LE approximation in one $r$-parametric line with $\sigma=10$ and $b=2.2$ ($6,000$ equidistant $r$-values in range $[0,300]$ are considered) parallel to those used to create the training, validation and test datasets. In top panel ${\rm{LE}}_1$ is shown, ${\rm{LE}}_2$ can be found in the middle panel, and finally, ${\rm{LE}}_3$ is in bottom panel. In each panel, the ground truth of the LEs (obtained with the algorithm in~\cite{WSSV85}) is in black, the mean of the predicted values with the $10$ networks is in red, and the uncertainty (that is, the interval $[${mean\,$\pm$\,standard deviation}$]$ obtained with the prediction of the $10$ networks) is in green. The obtained value of the Huber loss is equal to $0.091\pm0.005$. As already explained, to obtain this value, the prediction is performed with $10$ trained networks with the same architecture and data, but with different initilization of trainable parameters (weights and biases of the ANN). As the mean loss value and the standard deviation are small, we can conclude that the prediction is good enough.

At a first glance of the results in Figure~\ref{1DNonRandom_LS}, it can be seen that the parts where DL prediction seems to fail the most are those shaded in orange and purple. The orange one corresponds to orbits that converge to equilibrium points. If we analyse our normalization rule of the time series we can deduce that this fail is expected: the time series of equilibrium points are normalized to a random value between $0$ and $1$, so the CNN cannot extract information from them to predict the LE value (it is remarkable that the network assigns almost a constant value to the LEs of all the equilibrium point time series, so it has detected such kind of behaviour). As the LEs of an equilibrium point do not provide useful information among its negative sign in the first exponent (which is correctly predicted), we consider that the DL prediction is successful. The purple part corresponds to a region where long transient chaotic dynamics occur~\cite{BS07,BS09,TEL2008245}, and therefore, as we consider a quite short time series as data, it is logical that the DL technique can assume chaotic dynamics for them in some cases.

In general, the CNN that was trained only with a small number of $x$-time series is able to predict the three LE values correctly with a small uncertainty, and failing only in expected regions (as already explained) where it would be necessary to use complementary techniques to obtain a good result.

\subsubsection{2D Analysis}
\label{ss2}

\begin{figure*}[htb]
\centerline{\includegraphics[width=1.\textwidth]{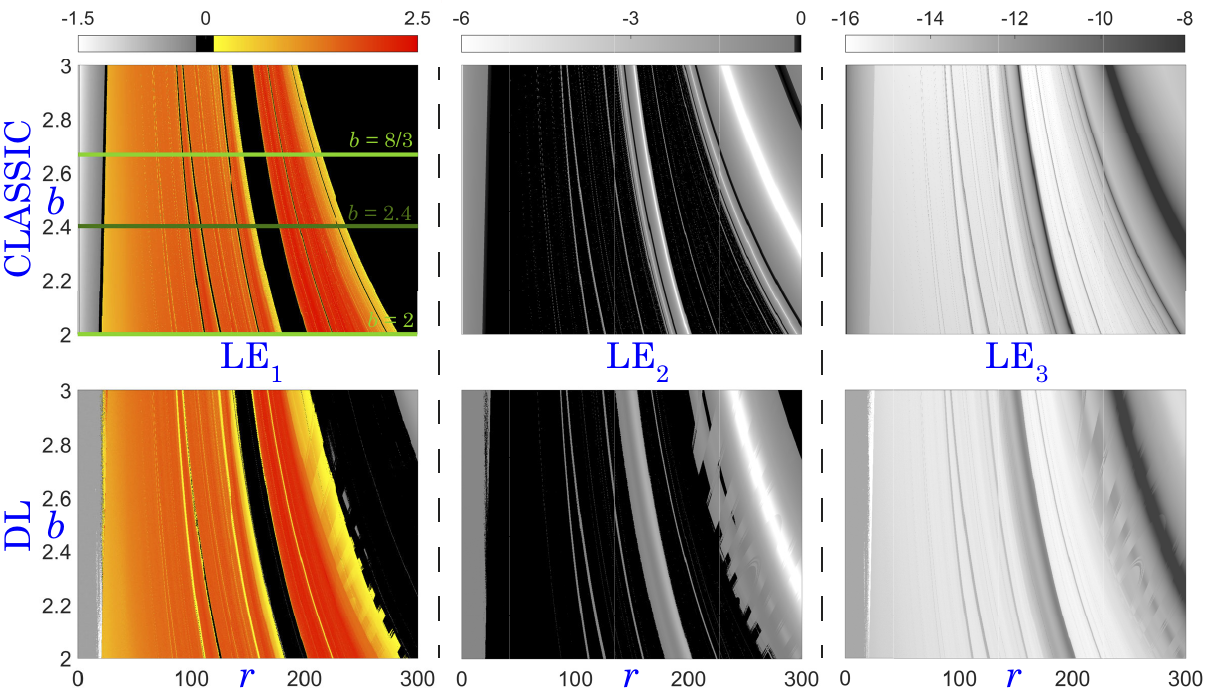}}
	\caption{2D biparametric analysis of Lyapunov Exponents in the Lorenz system ($\sigma=10$) when training with non-random data (Huber loss value $0.115\pm0.005$). From left to right, ${\rm{LE}}_1$, ${\rm{LE}}_2$ and ${\rm{LE}}_3$. From top to bottom, results with classical techniques and with DL techniques. Lines in the top-left panel correspond to lines from where training data (light green) and validation dataset (dark green) are obtained. (See the text for more details.)}
	\label{2DNonRandom_LS}
\end{figure*}

In Figure~\ref{2DNonRandom_LS} the CNN trained at the beginning of Section~\ref{sec:LS_NonRandom} with non-random data, based on three $r$-parametric lines (two for training and one for validation), has been used for LE prediction in the $(r,b)$-parametric plane ($1,000$ point values for each parameter, what makes a total of $10^6$ points). To obtain the time series for this biparametric analysis with the CNN, a transient integration is performed for $1,000$ unit times and later the integration is continued for $100$ more unit times (time step $0.01$ for all the integration). The input time series of length $1,000$ for the CNN are built with $1$ out of every $10$ of the last integration points. Remember that for the classical technique of LEs, transient integration is performed for $100,000$ time units (with time step $0.01$) and $10,001$ more time units (with time step $0.001$) are used to compute the LEs. As LEs are defined as a limit, with this classical technique, so much time is necessary to ensure a good LEs approximation. Therefore, with DL, less integration time is needed to approximate the full LEs spectrum.

In Figure~\ref{2DNonRandom_LS}, from left to right, ${\rm{LE}}_1$, ${\rm{LE}}_2$ and ${\rm{LE}}_3$. In the first row the results obtained with the classical technique in~\cite{WSSV85} are represented and second row corresponds to the DL prediction. To obtain such figure, black color is assigned to LEs with value around $0$, gray scale is used for negative values (different gray scales have been used in the color bars for a better interpretation: in all panels white is assigned to the smaller bar value and an enough dark gray is assigned to the largest negative value), and warm color gradation is used for positive values. Moreover, a minimum and a maximum value are fixed for each LE (that means that for example for the case of ${\rm{LE}}_1$, LEs with magnitude greater than $2.5$ are represented with the same color as Lyapunov Exponents of magnitude $2.5$, those with magnitude smaller than $-1.5$ are in the same color as those of magnitude $-1.5$, and the values in between follow the aforementioned gradation). The election of such minimum and maximum values does not affect the results as it is a usual way of representation for LEs.

At first sight, comparing the graphic results in Figure~\ref{2DNonRandom_LS} obtained with classical and DL techniques, it can be seen that, even predicting only with the short time series of the first variable $x$ of the Lorenz system, DL is able to reproduce quite well the magnitude of all the LEs. For the three LEs, the regions where the network seems to fail the most are the right boundary of the right chaotic region, and the upper right corner. In the aforementioned boundary, long transient chaos occurs. In spite of these small areas with non-precise approximations due to short time dynamics of the time series, the DL predictions would allow to perform a first dynamical analysis of the represented biparametric plane providing useful qualitative and quantitative approximations of the LEs. The Huber loss value is now $0.115\pm0.005$, not a large value considering the demanding task.

\begin{figure*}[htb]
\centerline{\includegraphics[width=1.\textwidth]{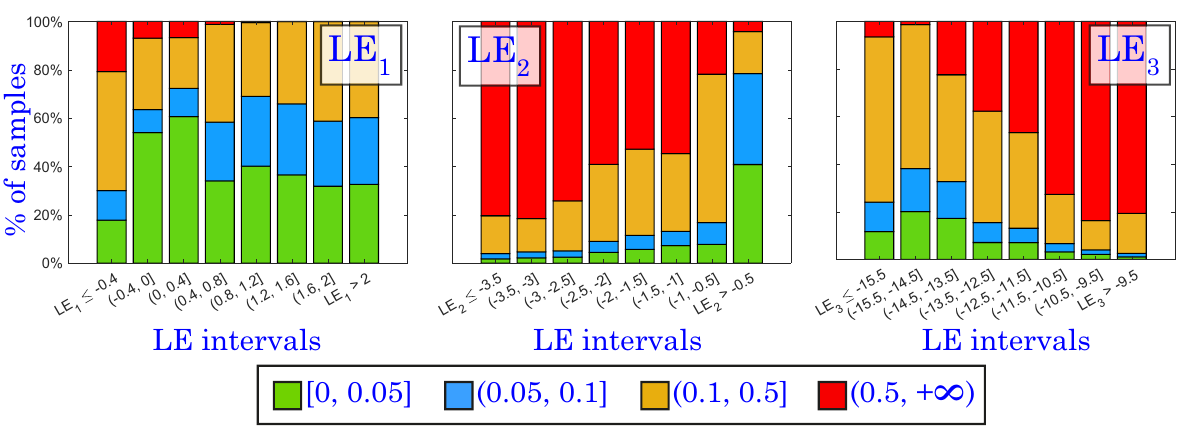}}
	\caption{Error analysis of Lyapunov Exponents prediction in an $(r,b)$-parametric plane of the Lorenz system (see Figure~\ref{2DNonRandom_LS}) when training with non-random data. From left to right, ${\rm{LE}}_1$, ${\rm{LE}}_2$ and ${\rm{LE}}_3$. Color code is given at the bottom. (See the text for more details.)} 
	\label{ErrorNonRandom_LS}
\end{figure*}

Figure~\ref{ErrorNonRandom_LS} shows the difference in absolute value between the predicted LE and that approximated by classical techniques~\cite{WSSV85} to quantitatively analyse the quality of the prediction. The error is given separately for each Lyapunov Exponent, from left to right, ${\rm{LE}}_1$ ,${\rm{LE}}_2$ and ${\rm{LE}}_3$. For each exponent, the samples are separated into different groups according to their LE value given by the classical technique (see the LE intervals in the horizontal axis of the figure). The differences in absolute value between the approximated LEs with the classical algorithm and with the CNN are computed. The percentage of samples (from the total number of samples on each LE interval) that belong to each of the error intervals ($[0, 0.05]$, $(0.05, 0.1]$, $(0.1, 0.5]$ and $(0.5, +\infty)$) is calculated and a color is assigned (green, blue, ochre and red, respectively for the mentioned error intervals, see the legend at the bottom of the figure) to obtain the error plots.

For ${\rm{LE}}_1$ (left panel), if the error is analysed, it can be seen that the largest errors are comitted for ${\rm{LE}}_1\leq -0.4$. As already mentioned in the 1D analysis, such LE values correspond to equilibrium points whose LE magnitude is not significant at all for a dynamical study. It is remarkable that for $\vert {\rm{LE}}_1\vert$ around $0$, more than half of the times, the error is less or equal than $0.05$. For ${\rm{LE}}_1>-0.4$ in almost all of the samples such prediction error is less or equal than $0.5$. Taking into account that the prediction is performed using one short time series of only one system variable, and without any extra dynamical information, predictions are quite accurate. For ${\rm{LE}}_2$ (middle panel), maybe the most important LE magnitude values for a dynamical analysis are those around $0$, which can provide some insight,  for instance, of period doubling bifurcations. The error analysis indicates that for $\vert {\rm{LE}}_2\vert$ around $0$, $40\%$ of the times the prediction is performed with an error not larger than $0.05$. Finally, the errors in ${\rm{LE}}_3$ prediction (right panel) are quite good since for all the LE intervals the error is less or equal than $0.5$ in more than $20\%$ of the times and this is a complicated approximation. We conclude that this DL technique allows us to predict all the LEs with only one variable short time series, and no other dynamical information of the system, with a good approximation (up to knowledge of the authors, using other techniques only the MLE is obtained when just one variable information is used).

\subsection{Random Data}\label{sec:LS_Random}

Last subsection was based on the supposition that we already have the Lyapunov exponents values in a few parametric lines and we used such data for training and validation. Now we suppose that we do not have any previous Lyapunov exponents data and so, we have to generate it to train the neural network. Therefore, we consider random data along the whole parametric plane we want to study in detail. For this approach, we randomly choose $24,000$ $(r,b)$-values ($\sigma = 10$) with $b\in[2,3]$ and $r\in[0, 300]$, and we obtain the $x$-time series (created as explained in Subsection~\ref{sec:LS_NonRandom}). Similar samples are deleted and the remaining ones are normalized as explained in Subsection~\ref{sec:LS_NonRandom}. Finally, $8,000$ samples are chosen randomly for training (batch size $128$), $2,000$ for validation (batch size $100$), and $2,000$ for test (batch size $100$). Notice that again only $8,000$ samples are directly related to training. The network architecture used for this LE prediction with random data is that explained in Subsection~\ref{sec:LS_NonRandom}.

The value of the Huber loss for training dataset is $0.054\pm0.003$. The corresponding values for validation and test datasets are $0.052\pm0.003$ and $0.057\pm0.003$, respectively. As the mean and standard deviation values are small enough, and as we consider that, however test mean value is larger than this for training, the difference is not remarkable to confirm that the DL technique suffers overfitting, it can be concluded that the training process has been successful. As indicated for the non-random case, the data used for the analyses with the trained CNN are just time series of $x$-variable of length $1,000$.

\subsubsection{1D Analysis}
\begin{figure}[htb]
\centerline{\includegraphics[width=0.5\textwidth]{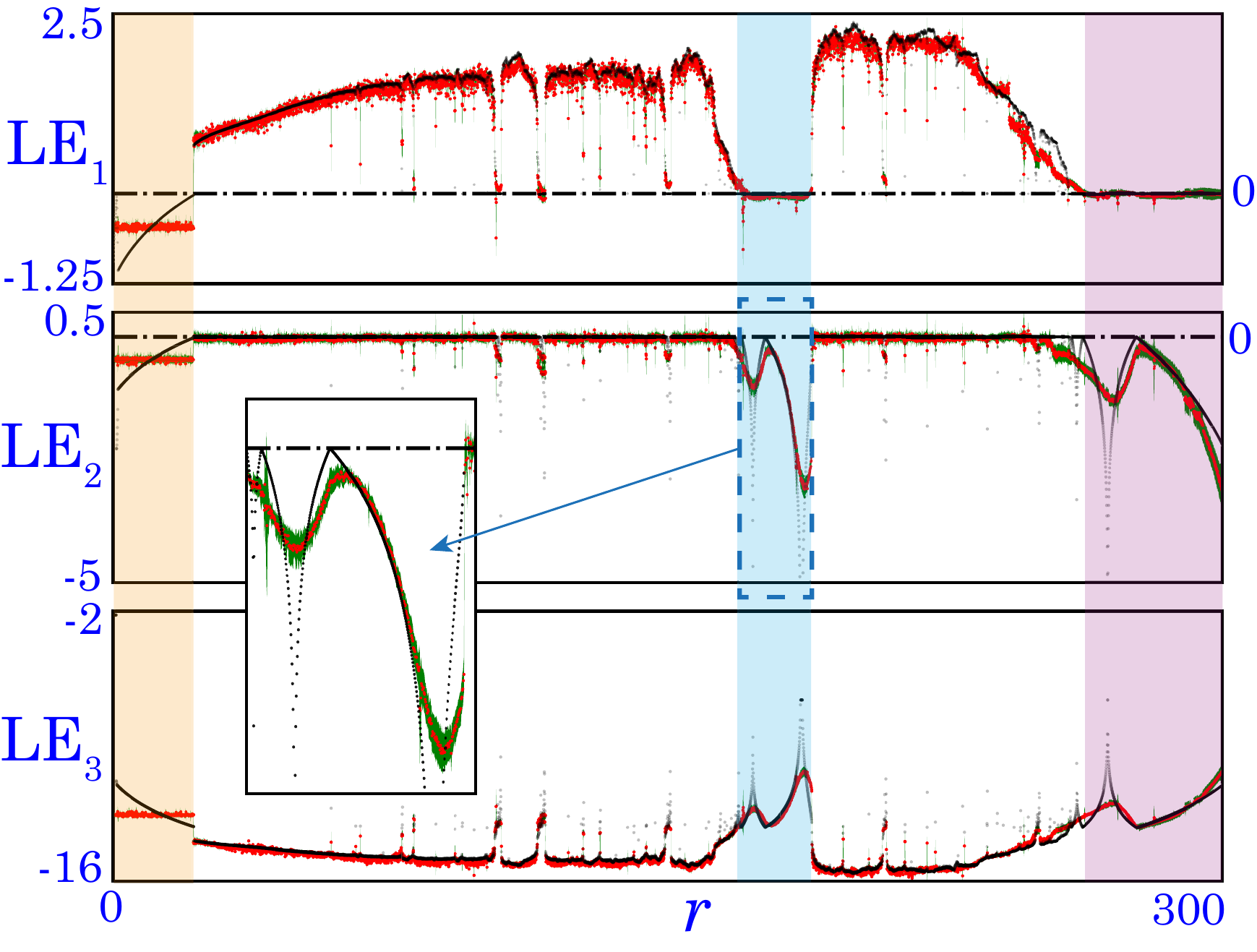}}
	\caption{1D parametric analysis ($\sigma=10$, $b=2.2$) of Lyapunov Exponents in the Lorenz system when training the CNN with random data (Huber loss value $0.055\pm0.003$). Orange and purple back colors correspond to regions where the DL technique seems to fail the most. Region shaded in blue is used to compare with previous analysis of Figure~\ref{1DNonRandom_LS}. (See the text for more details.)}
	\label{1DRandom_LS}
\end{figure}

As a first test, we consider the study of a one-parameter line.
In Figure~\ref{1DRandom_LS} the CNN trained with random data from the $(r,b)$-parametric plane has been used for LE approximation in one $r$-parametric line (with $6,000$ points) of such plane ($b=2.2$). As in Figure \ref{1DNonRandom_LS} (non-random case), in top panel ${\rm{LE}}_1$ is shown, ${\rm{LE}}_2$ is in middle panel, and finally, in bottom panel, ${\rm{LE}}_3$ is studied. The value of the Huber loss is $0.055\pm0.003$ (remember that numerical results are given as mean $\pm$ standard deviation for $10$ randomly initialized CNNs with the same architecture). The mean and standard deviation values are small, so we can confirm that the prediction task was successful.

As in the non-random data case (see Subsection \ref{sec:LS_NonRandom}), the orange part (equilibrium points) and the purple one (transient chaos) are those where the network seems to fail the most. However, if we compare the results obtained with each data creation technique, it can be seen that training with random data can obtain more accurate results in the critical region shaded in purple: in ${\rm{LE}}_1$ the uncertainty is smaller, and in ${\rm{LE}}_2$ and ${\rm{LE}}_3$, the shape followed by the DL results is more similar to that of the ground truth. This is an expected fact as a random sweeping to create the training data allows to obtain a bigger variability in dynamical behaviour for training than to restrict to just a few number of $r$-parametric lines.

Outside these problematic regions, the LE predictions of the CNN trained with random data are accurate and with small uncertainty. Let us compare the results in these regions with those of Subsection~\ref{sec:LS_NonRandom}. If the graphic representations are taken into account, the predictions of ${\rm{LE}}_1$ and ${\rm{LE}}_3$ seem to not present big differences. However, in the prediction of ${\rm{LE}}_2$ it is remarkable that the random dataset allows the network to predict in a more accurate way and with less uncertainty the LEs in the region shaded in blue (a magnification of this zone is shown in the figure). According to the loss function value, the predictions of the random case can be considered better as the given interval $[0.052,\,0.058]$ is closer to the origin than that of the non-random case $[0.086,\,0.096]$.

\subsubsection{2D Analysis}
\label{sec:r2}
\begin{figure*}[htb]
	\includegraphics[width=1.\textwidth]{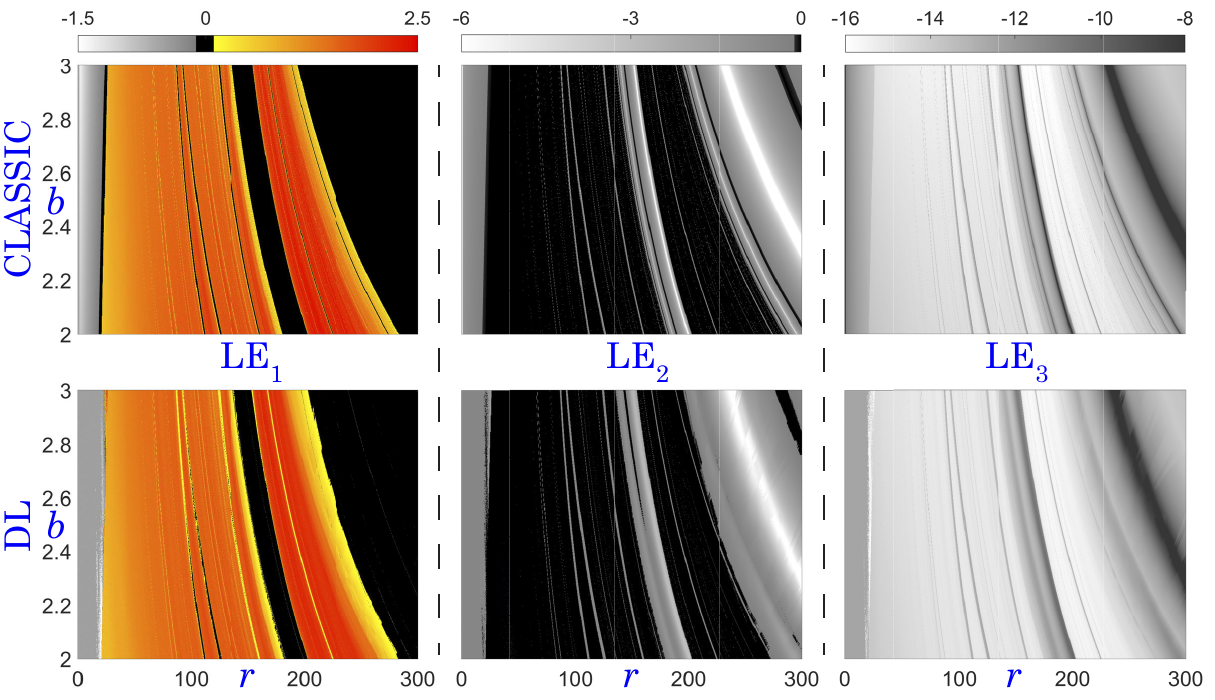}
	\caption{2D biparametric analysis of Lyapunov Exponents in the Lorenz system ($\sigma=10$) when training with random data (Huber loss value $0.079\pm0.006$). From left to right, ${\rm{LE}}_1$, ${\rm{LE}}_2$ and ${\rm{LE}}_3$. From top to bottom, results with classical techniques and with DL techniques. (See the text for more details.)}
	\label{2DRandom_LS}
\end{figure*}

Now we show the results when applying the trained networks for LE prediction in a biparametric plane ($\sigma=10$). Figure~\ref{2DRandom_LS} is equivalent to Figure~\ref{2DNonRandom_LS}, but now the CNN trained with random data has been used. Notice that only $8,000$ samples were taken for training (and $2,000$ for validation) from the $(r,b)$-biparametric plane in the figure. As for the non-random case, to compute the time series used by the DL technique, a transient integration is performed for $1,000$ time units and later the integration is continued for $100$ more time units (time step $0.01$ for all the integration). The input time series of length $1,000$ of the CNN are built with $1$ out of every $10$ of the last integration points. Notice that, with respect to classical technique of LEs, less time is needed to obtain the time series used by the CNN to compute the full LEs spectrum.

In Figure~\ref{2DRandom_LS} we compare the results obtained with classical and DL techniques, and it can be seen that, even predicting only with the first variable $x$ of the Lorenz system, DL is able to reproduce quite well the LE study of this region. At first sight, only the right boundary of the right chaotic region seems to present possible errors as it is blurred. This is the zone where transient chaos occurs and DL is expected to fail. However, if we compare it with the corresponding boundary in the non-random case (see Figure~\ref{2DNonRandom_LS}), it can be seen that using a random dataset the quality of the DL prediction is better for all the LEs. In fact, this improvement occurs in all the parametric plane in general. For example, the upper right corner whose values where predicted incorrectly in the non-random case (see Figure~\ref{2DNonRandom_LS}), now it is correctly represented. With a deeper visual analysis, the reader can realise that in the ${\rm{LE}}_2$ approximations, the predictions of the random case can help to detect some bifurcation lines (or dynamical changes) that the non-random technique did not allow (or not so clearly). For instance, the black line around $r$ parameter values between $150$ and $200$ (in the middle of the two big chaotic regions) in the ${\rm{LE}}_2$ panels appears in the DL panel of Figure~\ref{2DRandom_LS} for large values of $b$ and there is a darker gray for the smaller ones. In Figure~\ref{2DNonRandom_LS} this change cannot be seen so clearly. Another example is the black line around $r$ parameter values between $250$ and $300$. In the non-random case, there are some points in black, but there are not darker points that give idea of a continuous line. However, for the random prediction, even when black is not almost present, a continuous darker gray line highlights it. The value of the Huber loss for this random case is equal to $0.079\pm0.006$. This interval $[0.073,\,0.085]$ is closer to zero than that of the non-random data creation case $[0.110,\,0.120]$, so results are more accurate for random case as already shown in the 1D analysis.

\begin{figure*}[htb]
	\includegraphics[width=1.\textwidth]{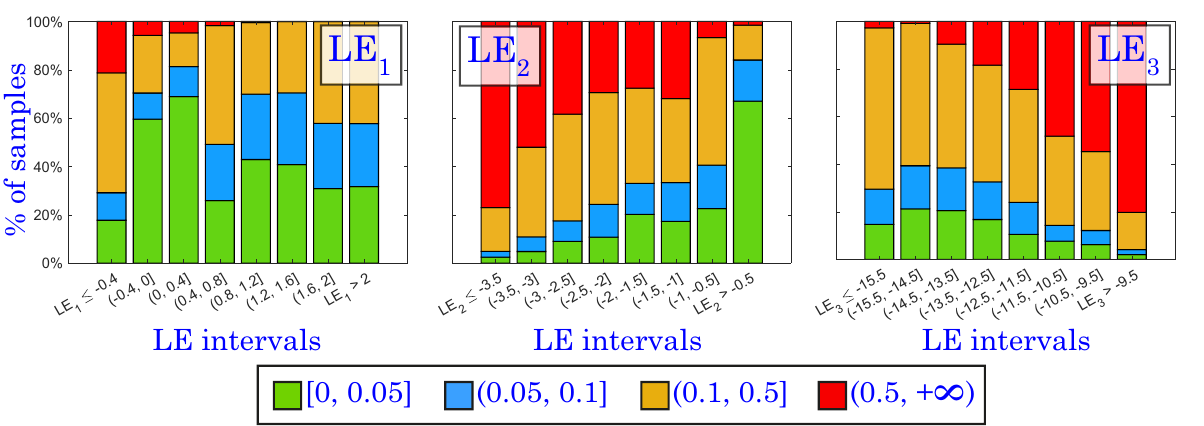}
	\caption{Error analysis of Lyapunov Exponents prediction in an $(r,b)$-parametric plane of the Lorenz system (see Figure~\ref{2DRandom_LS}) when training with random data. From left to right, ${\rm{LE}}_1$, ${\rm{LE}}_2$ and ${\rm{LE}}_3$. Color code is given at the bottom. (See the text for more details.)}
	\label{ErrorRandom_LS}
\end{figure*}

In Figure~\ref{ErrorRandom_LS} an error analysis equivalent to that of Figure~\ref{ErrorNonRandom_LS} is given for the random case. For ${\rm{LE}}_1$ (left panel), as in the non-random case, the biggest errors are committed for ${\rm{LE}}_1\leq -0.4$. It is remarkable that for $\vert {\rm{LE}}_1\vert$ around $0$, around $60\%$ of the times, the error is less or equal than $0.05$ (little improvement over the non-random case). For ${\rm{LE}}_1>-0.4$ the error is less or equal than $0.5$ for almost $100\%$ of the times. For ${\rm{LE}}_2$ (middle panel) and ${\rm{LE}}_3$ (right panel) the advantage of training with random data instead of non-random is remarkable. For ${\rm{LE}}_2$, when values are around $0$, the percentage of samples with an error less or equal than $0.05$ goes from $40\%$ in non-random case to more than $60\%$ in the current case. For ${\rm{LE}}_3$, the results are better in the random case. Notice that these predictions are quite good taking into account that the network had not much information for training and, up to the knowledge of the authors, there is not any other technique able to approximate ${\rm{LE}}_3$ in these conditions.

With all the study performed in the Lorenz system, it can be concluded that a good LE analysis can be performed with DL whether non-random or random data is used for training (with the second one providing better results). It is important to highlight that a small number of short time series are used to train (only $8,000$ for training, and $2,000$ more for validation, of length $1,000$), and only the $x$-variable of the system is used, which makes this a simple but powerful technique. Moreover, it is also a fast technique. As indicated in the part \textit{Lorenz system} of Table~\ref{timeTable}, it takes less than $1$ hour and $40$ minutes to compute a biparametric analysis from scratch with DL in the Lorenz system. Around $36$ minutes ($36\%$ of the total time used by the DL process) are needed to obtain the raw data that will be used to create train, validation and test datasets (CPU with parallel computing). Less than $40$ minutes ($40\%$ of the total DL time) are devoted to data selection (CPU), that is, to prepare such raw data and create the three mentioned datasets. To train one CNN (CUDA with PyTorch) less than $10$ minutes ($10\%$ of total time of DL process) are used (the results of the paper are obtained from $10$ random initialized CNNs, but as the standard deviation of the error is small, and as in the $1$D case the uncertainty seems to be small, it is expected that to use a unique trained CNN is enough to obtain good results). Finally, around $14$ minutes are used to obtain the full LEs spectrum with the trained CNN in a biparametric plane with dimension $1,000\times1,000$: only around $3$ seconds are devoted to the network prediction performed in CUDA with PyTorch, the remaining time is used to obtain the time series used as input to the network (CPU with parallel computing for some computations). Notice that most of the time used by DL ($76\%$) is focused on obtaining suitable data to train the network properly. With classical techniques (CPU with parallel computing), around $25$ hours are needed to perform such biparametric analysis. Therefore, comparing both techniques (classical and DL), with DL, time is reduced by $93.333\%$ approximately. In fact, once the network has been trained, time needed to obtain the biparametric analysis is less than $1\%$ of the time used by classical techniques. Moreover, if the time series are given and only the LEs are computed with the CNN, just $3$ seconds are needed to obtain the full LEs spectrum of the system.

\begin{table*}[htb]
	\begin{tabular}{rrrrrrr}
		\multicolumn{1}{l}{{\large{\textbf{\texttt{DEEP LEARNING}}}}} & \multicolumn{3}{c}{\textbf{Lorenz system}} & \multicolumn{3}{c}{\textbf{Coupled Lorenz system}} \\
		\hline
		 & \multicolumn{1}{r}{\textbf{Time}} & \multicolumn{1}{r}{\textbf{\% w.r.t. DL}} & \multicolumn{1}{r}{\textbf{\% w.r.t classical}} & \multicolumn{1}{r}{\textbf{Time}} & \multicolumn{1}{r}{\textbf{\% w.r.t. DL}} & \multicolumn{1}{r}{\textbf{\% w.r.t classical}}\\
		\hline
		\hline
		\textbf{Creation of raw data} & $36$ min & 36 \% & \multicolumn{1}{c}{-} & $68$ min & 49.275 \% & \multicolumn{1}{c}{-} \\
		\hline
		\textbf{Data selection} & $40$ min & 40 \% & \multicolumn{1}{c}{-} & $44$ min & 31.884 \% & \multicolumn{1}{c}{-} \\
		\hline
		\textbf{Training one CNN} & $10$ min & 10 \% & \multicolumn{1}{c}{-} & $10$ min & 7.246 \% & \multicolumn{1}{c}{-} \\
		\hline
		\hline
		\textbf{Biparametric analysis. Data} & $14$ min & 14 \% & \textbf{0.933 \%} & $16$ min & 11.594 \% & \textbf{0.494 \%} \\
		\hline
		\textbf{Biparametric analysis. Prediction} & $3$ s & 0.05 \% & \textbf{0.003 \%} &  $3$ s & 0.036 \% & \textbf{0.002 \%} \\
		\hline
		\hline
		\vspace{0.6cm}\textbf{Total time} & $1$ h $40$ min & 100 \% & \textbf{6.667 \%} & $2$ h $18$ min & 100\% & \textbf{4.259 \%} \\

		\multicolumn{1}{l}{{\large{\textbf{\texttt{CLASSICAL TECHNIQUE LEs}}}}} & \multicolumn{3}{c}{\textbf{Lorenz system}} & \multicolumn{3}{c}{\textbf{Coupled Lorenz system}} \\
		\hline
		 & \multicolumn{1}{r}{\textbf{Time}} & \multicolumn{1}{r}{\textbf{\% w.r.t. DL}} & \multicolumn{1}{r}{\textbf{\% w.r.t classical}} & \multicolumn{1}{r}{\textbf{Time}} & \multicolumn{1}{r}{\textbf{\% w.r.t. DL}} & \multicolumn{1}{r}{\textbf{\% w.r.t classical}}\\
		\hline
		\hline
		\textbf{Biparametric analysis. Whole process} & $25$ h & \multicolumn{1}{c}{-} & {100 \%} & $54$ h & \multicolumn{1}{c}{-} & {100 \%}\\
	\end{tabular}
	\caption{Time analysis for an LE biparametric study with DL and classical techniques. Top: Approximated time needed to perform a biparametric analysis with DL from scracth for the classical Lorenz system and a coupled Lorenz system. For each system, left column corresponds to time needed for each DL task (total time is given in the last row), middle column is for the percentage of time involved on each DL task, and right column is devoted to show the percentage of time used by DL respect to the time needed by classical technique of LEs for the same analysis. Bottom: Table with approximated time used by the classical technique of Lyapunov Exponents for both systems and the same biparametric analysis. Same meaning for columns as explained for DL. (See the text for more details.)}
	\label{timeTable}
\end{table*}

\section{LEs Approximation for the Coupled Lorenz System}\label{sec:coupledlorenz}

In this section, the two approaches used for the Lorenz system are applied to the system obtained coupling two almost identical Lorenz systems (see Equation~\ref{eqcls}) to approximate the full LEs spectrum with DL. Remember that on the first approach (non-random case) a CNN is trained (and validated) using time series from two (one) $r$-parametric lines, and on the second one (random case), the same architecture is trained (and validated) from scratch using information with random parametric values from an $(r,b)$-parametric plane ($\sigma=10$). As already mentioned, these two approaches will show that DL techniques can be used to expand a partial classical study of the system (first approach) or to do the analysis from random data (second approach). The performance of both approaches in the LE approximation in the coupled Lorenz system will be compared. The DL technique will allow us to locate hyperchaotic behaviour using only one variable data, in this case the $x_1$ variable of the first Lorenz system.

Again, the numerical results of the performance of the network are given as {mean\,$\pm$\,standard deviation} for $10$ randomly initialized networks.

As the study that we perform in the coupled Lorenz system is quite similar to the one performed in the Lorenz system in Section~\ref{sec:lorenz}, for simplicity, we present together the results obtained from both approaches.

\subsection{Data}

The creation of the training, validation and test dataset for the non-random case in the coupled Lorenz system is the same used for the Lorenz system. For the coupled system, parameters $k_1$ and $k_2$ are set to $0.1$.
From the samples satisfying $b\in\{2,8/3\}$, $8,000$ are chosen randomly for the training set (batch size $128$). For the data on the line with $b=2.4$, we select $2,000$ random points for validation (batch size $100$). Finally, $2,000$ random samples from the set $b=2.8$ are used for the test set (batch size $100$). The light green lines of the top-left panel of Figure~\ref{Plane_CLS_NR} correspond to the $r$-parametric lines used to obtain training data and the dark green line in the same plot corresponds to the one used to obtain validation data.

For the random case in the coupled Lorenz system, the creation of the training, validation and test dataset is the same used for the Lorenz system. From the random samples obtained in the parametric zone we want to study (using the same ranges that in Lorenz sytem study), $8,000$ samples are chosen randomly for training (batch size $128$), $2,000$ for validation (batch size $100$) and $2,000$ for test (batch size $100$). As in the Lorenz system study, $\sigma$ is set to $10$ and the samples are in the plane determined by $b\in[2,3]$ and $r\in[0,300]$. The parameters $k_1$ and $k_2$ of the coupled Lorenz system are set (both of them) to $0.1$.

The CNN architecture used for the prediction of the LEs spectrum for this dynamical system is that explained in Subsection~\ref{sec:LS_NonRandom}, just changing the last layer that now will have $6$ neurons. The coupled Lorenz system is a six dimensional system, so now the complete spectrum consists of six Lyapunov Exponents $\{{\rm{LE}}_1, \ldots, {\rm{LE}}_6\}$. Notice that again the data used for the analyses with the trained network are just time series of one variable (in this case $x_1$) of length $1,000$. Comparing with the classical algorithm that allows to obtain the whole LEs spectrum, with DL only one variable instead of six is used.

\subsection{1D Analysis}\label{sec:CLS_1D}
\begin{figure*}[htb]
	\centerline{\includegraphics[width=1.\textwidth]{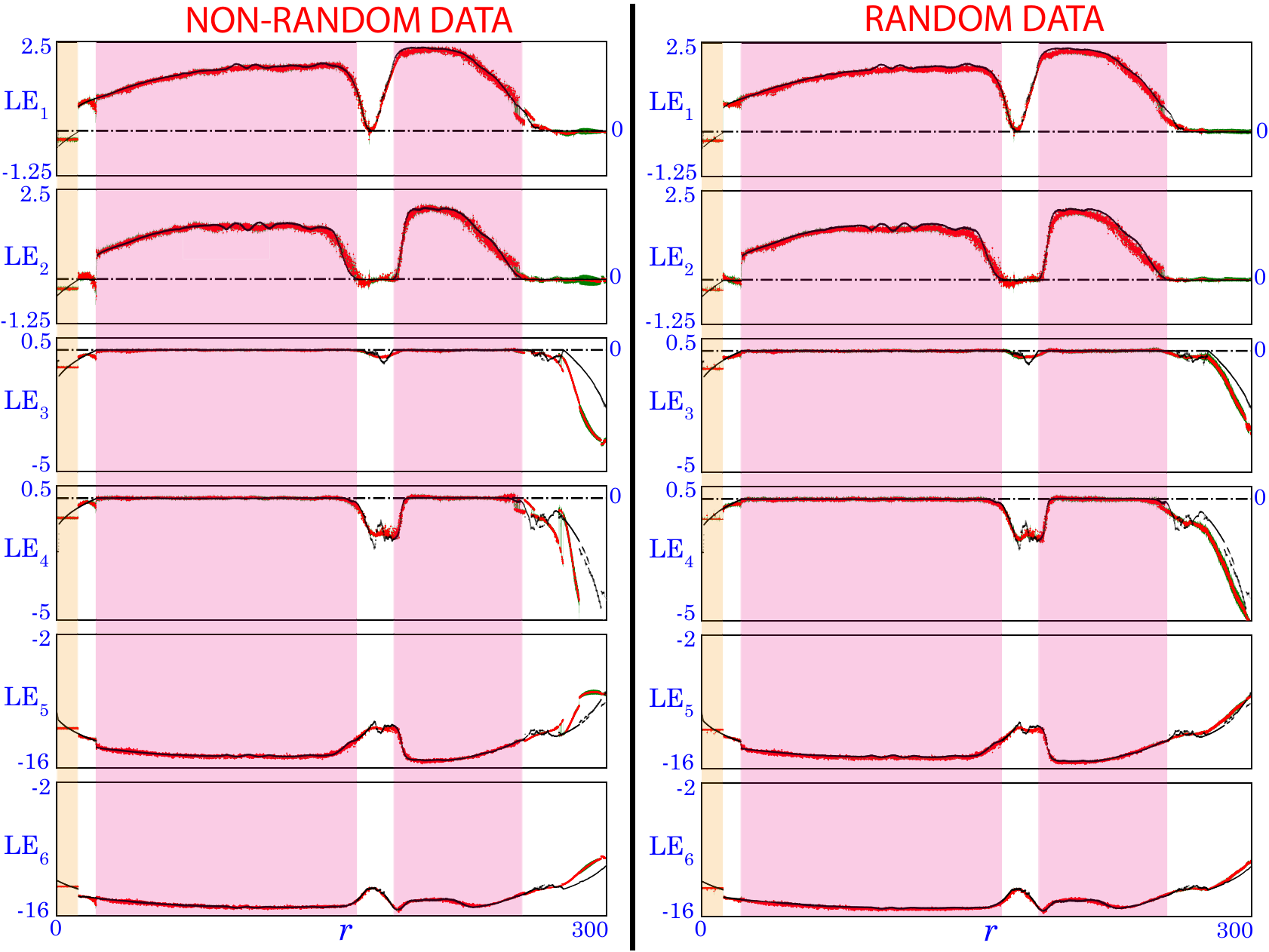}}
	\caption{1D analysis ($\sigma = 10$, $b=2.2$) of Lyapunov Exponents in the coupled Lorenz system when training the CNN with random and non-random data (Huber loss value in non-random approach is $0.274\pm 0.023$, and in random approach is $0.273\pm 0.027$). Region shaded in orange correponds to equilibirutm points where the network is expected to fail. Parts with pink back color are hyperchaotic regions. (See the text for more details.)}
	\label{Line_CLS}
\end{figure*}

In this subsection we make the analysis of one-parameter line ($b=2.2$ and $\sigma = 10$). In all panels of Figure~\ref{Line_CLS}, the ground truth of the LEs (obtained with the algorithm in~\cite{WSSV85}) is in black, the mean of the predicted values with the $10$ networks is in red, and the uncertainty (that is, the interval $[${mean\,$\pm$\,standard deviation}$]$ obtained with the predictions of the $10$ networks) is in green. In the figure we present the results for the complete spectrum with both approaches.

In the left column in Figure~\ref{Line_CLS}, the trained CNN has been used for LE approximation in one $r$-parametric line (with $6,000$ points) parallel to those used to create the datasets since we are in the non-random case. In the right column of the same figure the CNN trained with random data from the $(r,b)$-parametric plane has been used for LE approximation in the same $r$-parametric line (with $6,000$ points) of such plane.

The value of the Huber loss is $0.274\pm0.023$ in the non-random case and $0.273\pm0.027$ in the random approach. As a first sight, we can observe how both DL approaches provide very good approximated values of all the six Lyapunov Exponents using just one short time series!

We have marked the hyperchaotic intervals (${\rm{LE}}_1,\,{\rm{LE}}_2>0$) with a pink back color in all the panels.
In this coupled system there are large hyperchaotic regions that the DL algorithm easily detects. The interval with convergence to an equilibria is marked in orange color.
Note that in this region the algorithm does not provide exact values, but it is not necessary as no extra information is required, just that all of the LEs are negative.

We observe that training with random data gives more accurate results in the right part of the line since the DL results for the last four Lyapunov Exponents are closer to the values and the form of the ground truth. As in the Lorenz system, the use of a random sweeping to create the training data allows to obtain a bigger variability in dynamical behaviour for training (and therefore better approximations of the LE spectrum), than restricting these data to just a few number of $r$-parametric lines.

In any case, the CNN that was trained \emph{only} with a small number of $x_1$-time series is able to predict the complete LEs spectrum, giving correct values  with just a small uncertainty, using non-random and random training data. This is a quite remarkable result, as just one short time series is enough to approximate the six Lyapunov Exponents using a correctly trained CNN.

\subsection{2D Analysis}\label{sec:CLS_2D}

\begin{figure*}[p]
	\includegraphics[width=1.\textwidth]{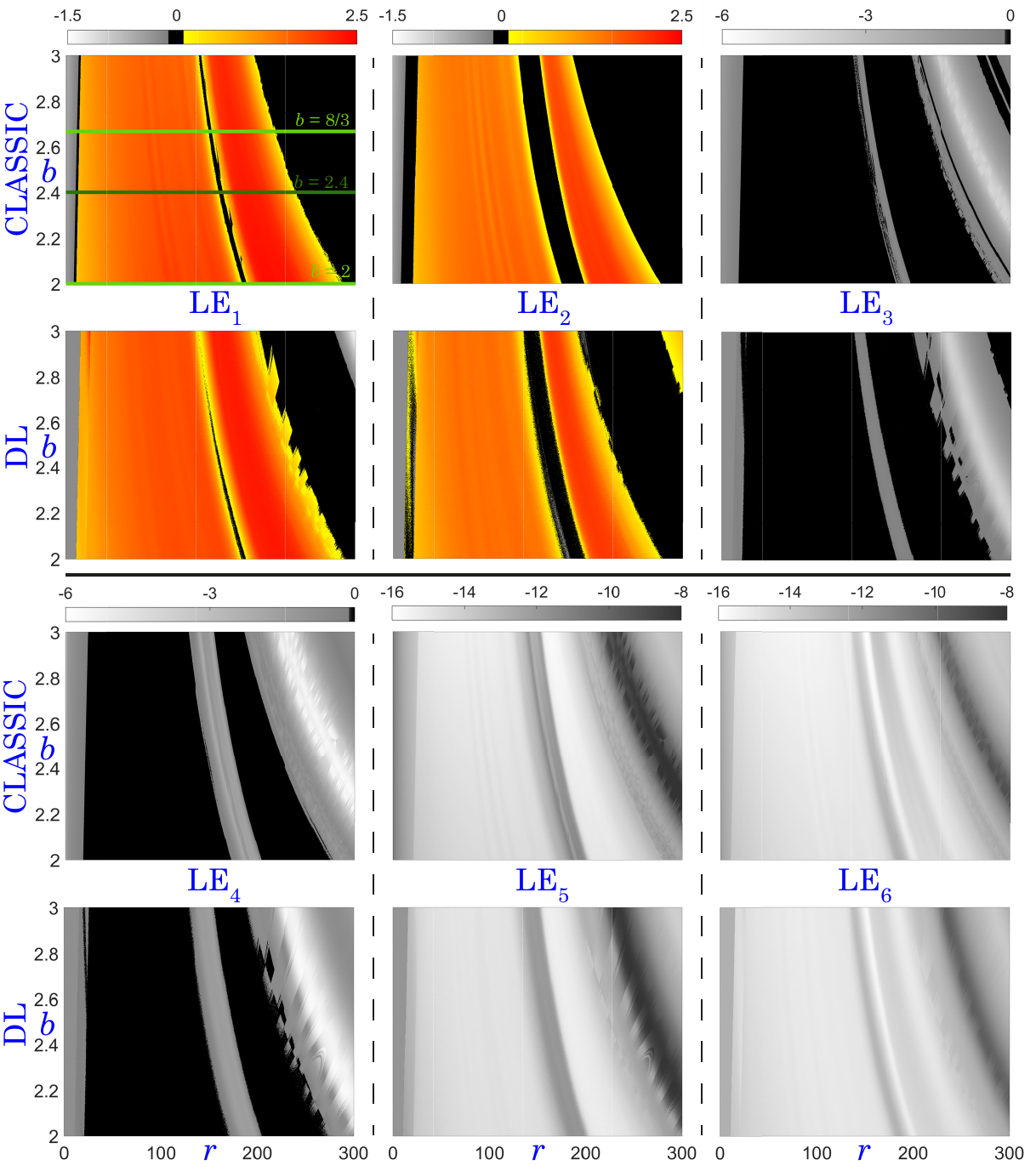}
	\caption{2D biparametric analysis of Lyapunov Exponents in the coupled Lorenz system ($\sigma=10$) when training with non-random data (Huber loss value is $0.046\pm 0.002$). (See the text for more details.)}
	\label{Plane_CLS_NR}
\end{figure*}

\begin{figure*}[p]
	\includegraphics[width=1.\textwidth]{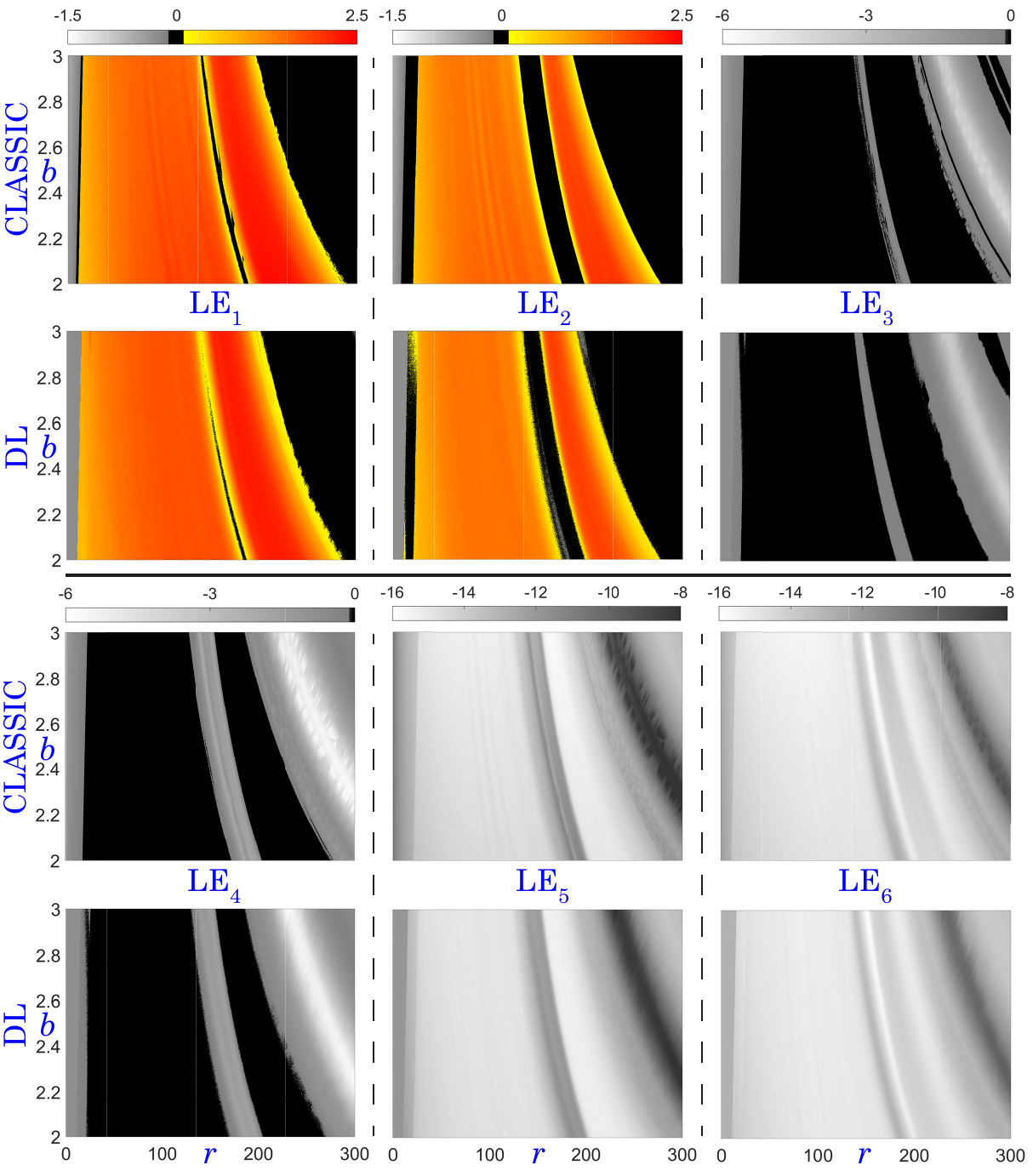}
	\caption{2D biparametric analysis of Lyapunov Exponents in the coupled Lorenz system ($\sigma=10$) when training with random data (Huber loss value is $0.023\pm 0.004$). (See the text for more details.}
	\label{Plane_CLS_R}
\end{figure*}

In this subsection we show the results in a biparametric plane ($\sigma = 10$) of the coupled Lorenz system. Figures~\ref{Plane_CLS_NR} and \ref{Plane_CLS_R} corresponding to the non-random and random approach, respectively, are equivalent to Figures~\ref{2DNonRandom_LS} and \ref{2DRandom_LS} of the Lorenz system. As for that system, to obtain the time series used by the DL technique, a transient integration is performed for $1,000$ time units and later the integration is continued for $100$ more time units (time step $0.01$ for all the integration). The input time series of the CNN of length $1,000$ are built with $1$ out of every $10$ of the last integration points. Therefore, with respect to classical techniques, less time is needed to obtain the time series that will be used to compute the LEs.


In both Figures~\ref{Plane_CLS_NR} and \ref{Plane_CLS_R} we observe how the DL techniques work correctly. 
The non-random approach has more problems in the transition areas from regular to chaotic behaviour and viceversa, but in any case the global result is quite good.

The random data approach, done as in Section~\ref{sec:r2}, provides a  much better approximation of the real LEs values giving a very similar result than the Lyapunov Exponents obtained from the standard algorithm.

In fact, not only the use of just short $x_1$-time series allows to detect chaotic behaviour, but also hyperchaotic one. And moreover, the technique provides a biparametric study of the regions with hyperchaoticity as shown in the pictures with the ${\rm{LE}}_2$ results, where in yellow-red scale the positive values of the second Lyapunov Exponent are represented.

In any case, the CNN that was trained only with a small number of $x_1$-time series is able to predict the complete LEs spectrum, giving correct values  with a small uncertainty.

\begin{figure*}[htb]
	\centerline{\includegraphics[width=0.95\textwidth]{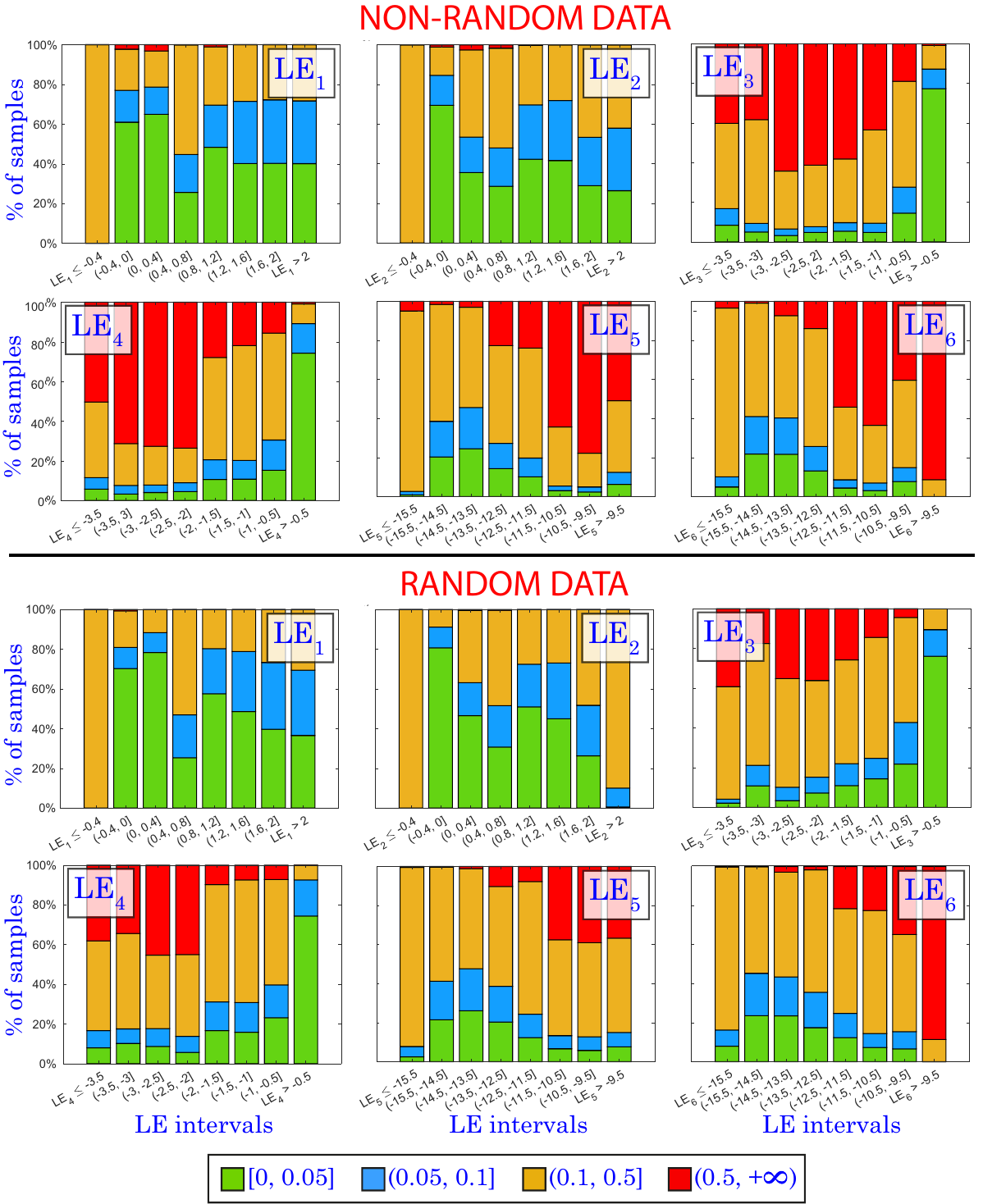}}
	\caption{Error analysis of Lyapunov Exponents prediction in an $(r,b)$-parametric plane of the coupled Lorenz system when training with non-random and random data. Color code is given at the bottom. (See the text for more details.)}
	\label{Error_CLS}
\end{figure*}

To complete the study, in Figure~\ref{Error_CLS} an error analysis is provided to analyse quantitatively the quality of the prediction taking into account the difference in absolute value between the predicted LE and the value approximated by classical techniques. This analysis is given separately for each Lyapunov Exponent. As in previous section, the samples are separated for each exponent into different groups according to their LE value given by the classic technique (see the LE intervals in the horizontal axis of the figure). For each LE interval, the percentage of samples that belong to each of the error intervals ($[0, 0.05]$, $(0.05, 0.1]$, $(0.1, 0.5]$ and $(0.5, +\infty)$) is calculated and a color is assigned for these error intervals to obtain the error plots (green, blue, ochre and red, respectively).
As expected, the advantage of training with random data instead of non-random is remarkable for the performance of the high order Lyapunov Exponents (the error intervals representing the biggest errors are less significative in the random approach).
Therefore, we can extract similar conclusions as in the case of the Lorenz system, and we conclude that this DL technique allows us to predict all the LEs with only one variable short time series (no other dynamical information of the system is needed), with a good approximation.

From the performed studies in the coupled Lorenz system, it can be concluded that a good LE analysis can be obtained with DL whether non-random or random data is used for training (with the second one providing better results as in the Lorenz system studies). It is remarkable that a small number of short time series have been used to train ($8,000$ samples for training, and $2,000$ more for validation, of length $1,000$), and just one of the six variables of the system is used. So, this is a simple but powerful technique. In addition, it is also a fast technique. As indicated in the part \textit{Coupled Lorenz system} of Table~\ref{timeTable}, to obtain a biparametric analysis from scratch with DL in this coupled Lorenz system, less than $2$ hours and $18$ minutes are needed. Around $68$ minutes ($49.275\%$ of the total DL time) are devoted to obtain the raw data used later to create train, validation and test datasets (CPU with parallel computing). Less than $44$ minutes ($31.884\%$ of the total time used by the DL process) are needed for data selection, that is, to prepare such raw data and create the three mentioned datasets (CPU). To train one CNN (CUDA with PyTorch) less than $10$ minutes ($7.246\%$ of total DL time) are used (the results of the paper are obtained from $10$ random initialized CNNs, but because of the obtained results for this coupled system, it is expected that to use a unique CNN is enough to obtain good results). Finally, around $16$ minutes are used to obtain the full LEs spectrum with the trained CNN in a biparametric plane with dimension $1000\times1000$: only $3$ seconds ($0.036\%$ respect total DL time) are needed by the network prediction in CUDA with PyTorch, the remaining time ($11.594\%$ of total DL time) is used to obtain the time series used as input to the network (CPU with parallel computing for some computations). With classical techniques (CPU with parallel computing), almost $54$ hours are needed to obtain such biparametric analysis, so comparing both techniques (classical one and DL), with DL, time is reduced by $96\%$ approximately. In fact, once the network has been trained, time needed to obtain the biparametric analysis is less than $0.5\%$ of the time used by classical techniques. Moreover, in the case in which the time series are given and only the LEs have to be computed with the trained CNN, just $3$ seconds are needed to obtain the full LEs spectrum of the system.

\section{Conclusions}\label{sec:Conclusions}
The Lyapunov Exponents spectrum of a dynamical system is arguably one of its most fundamental properties as it permits to characterize the dynamics of the system. Its computation can be highly computationally expensive, specially if one focuses on a classification problem in a parameter plane. But this information can provide a global panorama on the dynamics, and so it is quite important.

In this paper, a well-known Deep Learning network (Convolutional Neural Network, CNN) has been built and trained to carry out the approximation of the complete Lyapunov Exponents spectrum of a dynamical system.
The training process is performed using as data the complete Lyapunov spectrum of a small number of points in the parametric space we work on, but once trained, the network only needs short time-series in just one variable of the system, which means a large reduction in time and memory for the approximation of the LEs. The methodology has been applied in two test problems: the Lorenz system and the coupled Lorenz system.

For the Lorenz system, we have used the trained networks to study the behaviour of an $r$-parametric line and a biparametric plane of the parameter space.
For the coupled Lorenz system, we study the behaviour on the same $r$-parametric line and biparametric plane as in the isolated Lorenz model, but now as the system has dimension six we use the network to approximate the six Lyapunov Exponents. We highlight that the training process uses just a few lines of one-parameter data or a short number of random points to create a network capable of performing biparametric studies. This is a remarkable result that shows us the power of DL techniques in dynamical systems studies.

For the biparametric study of the Lorenz system, around $25$ hours are needed to perform such analysis with classical techniques, while with the CNN less than $2$ hours are needed for the same task, that is, a $93.333\%$ time is saved with DL techniques (if time series are given, the time of the predictions is just $3$ seconds). In the case of the biparametric study of the coupled Lorenz system, almost $54$ hours are needed to obtain this analysis with classical techniques, while with the CNN just over two hours are necessary, therefore, a $96\%$ time is saved with DL (if time series are given, the prediction time is just $3$ seconds).

We conclude that Deep Learning can be used to not only analyse the behaviour (regular, chaotic or hyperchaotic) of a dynamical system, but also to quantify the values of the Lyapunov Exponents spectrum, that is, to go further to a classification problem. Our results show that even dense $2$D parametric studies can be carried out in a very reasonable time using data from just a small portion of the global phase space. However, a deeper study would be necessary to know how far we can go using these techniques in this and other dynamical systems tasks. In summary, we have demonstrated that inference of the full Lyapunov Exponents spectrum from a short single variable time series is possible with DL and it is robust.

\section*{Acknowledgments}
RB, AL, AMC and CMC have been supported by the Spanish Research projects PID2021-122961NB-I00. RB, AMC and CMC have been supported by the European Regional Development Fund and Diputaci\'on General de Arag\'on (E24-23R). RB has been supported by the European Regional Development Fund and Diputaci\'on General de Arag\'on (LMP94-21). AL has been supported by the European Regional Development Fund and Diputaci\'on General de Arag\'on (E22-23R).

\section*{CRediT Authorship Contribution Statement}
C. Mayora-Cebollero: Conceptualization, Formal analysis, Investigation, Methodology, Validation, Visualization, Writing–original draft, Writing–review \& editing. A. Mayora-Cebollero: Conceptualization, Investigation, Supervision, Writing–review \& editing. \'A. Lozano: Investigation, Supervision, Writing–review \& editing.
R. Barrio: Conceptualization, Formal analysis, Project administration, Methodology, Supervision, Writing–review \& editing.

\section*{Declaration of Competing Interest}
 There is no conflict of interest between the authors and other persons organizations.

\section*{Data Availability}
Data available on request from the authors.

\bibliographystyle{elsarticle-num}

\end{document}